\documentclass[conference]{IEEEtran}

\usepackage{utils/packages}
\usepackage{textcomp}

\def\BibTeX{{\rm B\kern-.05em{\sc i\kern-.025em b}\kern-.08em
    T\kern-.1667em\lower.7ex\hbox{E}\kern-.125emX}}

\newtheorem{theorem}{Theorem}

\newtheorem{assumption}{Assumption}
\newtheorem{remark}{Remark}

\DeclareMathOperator*{\argmin}{argmin}

\algnewcommand{\Initialize}[1]{%
  \State \textbf{Initialize:}
  \Statex \hspace*{\algorithmicindent}\parbox[t]{.8\linewidth}{\raggedright #1}
}

\begin{document}

\title{Network-GIANT: Fully distributed Newton-type optimization via harmonic Hessian consensus}


\author{\IEEEauthorblockN{Alessio Maritan$^*$, Ganesh Sharma$^\dag$, Luca Schenato$^*$, and Subhrakanti Dey$^\ddag$}
\IEEEauthorblockA{$^*$Department of Information Engineering, University of Padova, Italy \\
$^\dag$Hamilton Institute, Maynooth University, Ireland \\
$^\ddag$Department of Electrical Engineeringe, Uppsala University, Sweden \\
Email: alessio.maritan@phd.unipd.it, ganesh.sharma.2019@mumail.ie, schenato@dei.unipd.it, subhrakanti.dey@angstrom.uu.se}
}
\maketitle
\thispagestyle{empty}


\begin{abstract}
This paper considers the problem of distributed multi-agent learning, where the global aim is to minimize a  sum of local objective (empirical loss) functions through local optimization and information exchange between neighboring nodes. We introduce a Newton-type fully distributed optimization algorithm, Network-GIANT, which is based on  GIANT, a Federated learning algorithm that relies on a centralized parameter server. The Network-GIANT algorithm is designed via a combination of gradient-tracking and a Newton-type iterative algorithm at each node with consensus-based averaging of local gradient and Newton updates.
The resulting algorithm is efficient in terms of both communication cost and run-time, making it suitable for wireless networks.
We prove that our algorithm guarantees semi-global and exponential convergence to the exact solution over the network assuming strongly convex and smooth loss functions. We provide empirical evidence of the superior convergence performance of Network-GIANT over other state-of-art distributed learning algorithms such as Network-DANE and Newton-Raphson Consensus.

\end{abstract}

\section{Introduction}
\IEEEPARstart{W}ith billions of IoT devices, one now can have access to huge volumes of data, which can also be high dimensional. It is difficult to process such large volumes of data on a single machine. Distributed learning and optimization algorithms are therefore attracting significant attention recently for performing machine learning tasks. Distributed machine learning (DML) involves training with local data sets across using multiple participating machines that work together in a coordinated manner to learning a global model, without sharing raw data. There are two main frameworks for DML: a) federated learning, which involves a parameter server acting as a coordinator in a master/worker setting \cite{McMahanFL}, and b) networked learning, which has no parameter server and nodes communicate with their single-hop neighbours in a peer-to-peer setting to reach a consensus over shared information. Our focus in this paper is on the latter, fully distributed networked learning problem.

\par {\bf Related Work}:
There has been extensive work in DML in the federated learning framework both with gradient-based optimization (see \cite{GENERAL_FL} for a survey) and Hessian-based optimization exploiting approximate Newton-type methods e.g. DANE \cite{shamir2014communication}, GIANT \cite{wang2018giant}, DONE\cite{dinh2022done}, DANDA \cite{sharma2022analog}, and FedNL \cite{safaryan2021fednl}. The networked setting has been explored with gradient-based optimization in a number of works, e.g. \cite{nedic2009distributed,jakovetic2014fast,yuan2016convergence,nedic2017achieving}, whereas networked version of the Alternating Direction Method of Multipliers (ADMM) have been explored in \cite{shi2014linear,chang2014multi}. 
In the context of Newton-type second-order methods, the Newton-Raphson method has been extended to the network setting in the works \cite{bof2018multiagent} and \cite{varagnolo2015newton}. In particular, the algorithm proposed in the latter either requires the communication of the entire Hessian matrices, which can be prohibitively expensive, or exploits only the diagonal of Hessian matrices, possibly leading to slow convergence for skewed objective functions.
The authors in \cite{mokhtari2016network} and \cite{bajovic2017newton} have proposed Network-Newton methods based on an approximate Newton step utilizing truncated Taylor series expansion of the Hessian inverse, but with a penalty-based modified cost function. In \cite{zhang2021newton} the authors have proposed a Newton tracking algorithm where each node updates its local variable along a local Newton direction modified with neighboring and historical information. The authors in \cite{li2020communication} develop a network version of the DANE algorithm called Network-DANE. The latter requires to solve an inner minimization problem a each iteration, which can be computationally demanding and can involve hyperparameter tuning for the sub-solver.
\par
\textbf{Main Contributions}: In this paper we propose a fully decentralized (networked) version of the GIANT \cite{wang2018giant} algorithm, termed as Network-GIANT, which operates without a central parameter server acting as of communication and information aggregator for all nodes. Instead, the learning takes place over a network (connected graph) of nodes where each node runs a two-steps procedure, leveraging gradient tracking and average consensus. Network-GIANT enjoys a low communication overhead of 2 $d$-dimensional vectors at each iteration. We formally prove the semi-global exponential convergence of Network-GIANT to the exact optimal solution. Finally, we present a set of empirical studies comparing Network-GIANT to other state-of-the-art network learning algorithms, demonstrating its superior convergence performance.
\par
\textbf{Organization}: 
The remainder of the  paper is organized as follows. Section 2 describes the problem and motivation behind our work. The proposed Network-GIANT algorithm is presented in Section 3, and Section 4 details the convergence analysis of the algorithm. Section 5 presents numerical simulations on two real datasets, followed by concluding remarks in Section 6.

\section{Problem Formulation}
Consider a set of agents $\mathcal{N} = \{1, \dots, n\}$ where each member $i$ has access to a possibly private local cost function $f_i: \mathbb{R}^d \to \mathbb{R}$. We are interested in solving in a distributed fashion the unconstrained optimization problem
\begin{equation}
    f(x^\star) = \min_{x \in \mathbb{R}^d} \left\{ f(x) = \frac{1}{n} \sum_{i=1}^n f_i(x) \right\}
    \label{eq:problem_formulation}
\end{equation}
where each $f_i(x)$ is strongly convex. This general formulation suits several problems that frequently arise in machine learning, in particular the task of empirical risk minimization. In the latter, $x \in \mathbb{R}^d$ parametrizes a model to be learnt and each $f_i(x)$ is a suitable loss function evaluated over the data samples possessed by agent $i$. 
In the context described above, we want to devise a second-order distributed optimization algorithm which exploits the convexity of the objective functions to achieve faster convergence. Our goal is to reach the exact solution, differently from \cite{bajovic2017newton} and \cite{mokhtari2016network} which can only guarantee convergence to a neighborhood of the latter.

\textit{Notation}: For a generic variable $x$, we use $x_i^k$ to denote the local copy of $x$ possessed by agent $i$ at iteration $k$. $\left\| \cdot \right\|$ is the Euclidean norm, $I$ is the identity matrix and the superscript $^T$ indicates the transpose of the argument. $\mathbb{1}$ and $\mathbb{0}$ are $d$-dimensional vectors whose components are all equal to 1 and 0, respectively.

\bigskip \noindent
\textbf{GIANT.} In the Federated Learning setting, the authors of \cite{wang2018giant} addressed a class of empirical risk minimization problems and proposed GIANT, a Newton-type method based on the update
\begin{equation}
    x^{k+1}_{\text{GIANT}} = x^k - (H^k)^{-1} \left( \frac{1}{n} \sum_{i=1}^n \nabla f_i(x^k) \right) ,
    \label{eq:GIANT_update_rule}
\end{equation}
\begin{equation*}
H^k = \left( \frac{1}{n} \sum_{i=1}^n \left[ \nabla^2 f_i(x^k)^{-1} \right] \right)^{-1} .
\end{equation*}
In comparison, the traditional Newton-Raphson recursion relies on the true Hessian of $f(x)$, given by the arithmetic mean of the local Hessians. According to \cite{wang2018giant}, if the information is spread-out among the data samples then the harmonic mean $H^k$ is very close to the arithmetic mean, and therefore \eqref{eq:GIANT_update_rule} becomes a good approximation of the Newton update. In practice, a line search procedure is needed to select a suitable step-size and prevent algorithmic divergence.

GIANT is designed for the master-worker architecture, where all the agents are connected in a star topology to a central aggregator node.
The job of the latter is to gather the local information and compute the averages $\nabla f(x^k)$ and $\eta$, which are respectively the mean of the local gradients and the mean of the approximate Newton directions $\eta_i = \nabla^2 f_i(x^k)^{-1} \nabla f(x^k)$ $i = 1, \dots, n$. The pseudocode refers to the version of GIANT in which the step-size is chosen by means of the following line search procedure. First, each worker $i$ evaluates its local function at $x_i^k - \alpha_j \eta$ for a fixed set of candidate step-sizes $\{ \alpha_j \}$. Then the master selects a step-size $\alpha^\star$ that satisfies the Armijo–Goldstein condition and the workers use it to update their local variable.


Overall, each iteration requires the transmission for each agent of 4 $d$-dimensional vectors, a vector whose size is the cardinality of $\{ \alpha_j \}$ and a scalar, resulting in a competitive communication complexity.

Despite being provably more efficient than several first-order and second-order methods, the applicability of GIANT is limited to the master-worker scenario. To overcome this shortcoming, in the next section we introduce Network-GIANT, which allows to implement the update formula \eqref{eq:GIANT_update_rule} in a general network of agents, without the need of central coordinator.

\renewcommand{\thealgorithm}{}
\setlength{\columnseprule}{0.4pt}
\begin{algorithm}
\caption{GIANT with line search \cite{wang2018giant}}\label{alg_GIANT}
\begin{algorithmic}
\Initialize{
    Choose $x^0 \in \mathbb{R}^d$, set $x_i^0 = x_0 \ \forall i = 1, \dots, n$. \\
    Choose the set of candidate step-sizes $\{ \alpha_j \}$.}
\For{each iteration $k = 0, 1, 2, \dots$}
   \State $\nabla f(x^k) = \frac{1}{n} \sum_{i=1}^n \nabla f_i(x^k)$ \Comment{at the master}
   \For{each worker $i = 1, \dots, n$}
       \State $\eta_i = \nabla^2 f_i(x^k)^{-1} \nabla f(x^k)$
   \EndFor
   \State $\eta = \frac{1}{n} \sum_{i=1}^n \eta_i$ \Comment{at the master}
   \For{each worker $i = 1, \dots, n$}
       \State $f_{i,\alpha_j} = f_i(x_i^k - \alpha_j \eta) \quad \forall \alpha_j \in \{ \alpha_j \}$
   \EndFor
   \State $\alpha^\star = \argmin_{\alpha_j}  \frac{1}{n} \sum_{i=1}^n f_{i,\alpha_j}$ \Comment{at the master}
   \For{each worker $i = 1, \dots, n$}
       \State  $x_i^{k+1} = x_i^k - \alpha^\star \eta$
   \EndFor
\EndFor
\end{algorithmic}
\end{algorithm}

\section{Network-GIANT}
The GIANT algorithm assumes the presence of a master node which can gather all the local information, compute the mean and send it back to the workers in one shot. In this section we show how to remove the dependency from such central entity, using the tools of distributed average consensus and gradient tracking. We first provide a formal definition of the multi-agent setting under consideration, and then we show how to modify GIANT to make it suitable for peer-to-peer distributed optimization.


We consider the communication network represented by the connected and time-invariant graph $\mathcal{G} = (\mathcal{N}, \mathcal{E})$, where the vertex set $\mathcal{N}$ is the set of the agents. The edges $\mathcal{E}$ are the available bi-directional communication links, and each node can only communicate with its single-hop neighbors.
The network is associated with a consensus matrix $P \in \mathbb{R}^{n \times n}$, whose generic entry $p_{ij}$ is positive if edge $(i, j) \in \mathcal{E}$ and zero otherwise. The mixing matrix $P$ is symmetric and doubly stochastic, i.e. such that $P \mathbb{1} = \mathbb{1}$, $\mathbb{1}^T P = \mathbb{1}^T$ and the null space of $I-P$ is span$(\mathbb{1})$. It is possible to build a matrix that satisfies these requirements without complete knowledge of the graph topology using e.g. the Metropolis weights \cite{xiao2007distributed}.

In this setup, the centralized averages originally computed by the master can be replaced by two average consensus blocks. Consensus is an iterative process where each step is composed of two stages: a communication phase, in which each agent $i$ transmits its local variable $x_i$ to its neighbors, followed by a local update, where $x_i$ is updated using a weighted average of the information received by node $i$. 
The weights are given by the consensus matrix $P$, and this process continues until all agents converge to the same value, which corresponds to the arithmetic mean of the initial values of the local variables.

However, standard average consensus is not suited to efficiently estimate the gradient of the global cost. We then resort to the technique of average tracking consensus, introducing the local variable
\begin{equation}
    w^{k+1}_i = \sum_{j=1}^n p_{ij} \left( w_j^k + \nabla f_j (x_j^k) - \nabla f_j (x_j^{k-1}) \right) .
    \label{eq:w_gradient_tracking}
\end{equation}
Thanks to the doubly stochasticity of the matrix $P$, choosing any initialization that satisfies $w_i^0 = \nabla f_i (x_i^{-1})$ $\forall i \in \mathcal{N}$ we get the useful property
\begin{equation}
    \frac{1}{n} \sum_{i=1}^n w^{k+1}_i = \frac{1}{n} \sum_{i=1}^n \nabla f_i (x_i^k) .
\end{equation}
As the consensus process proceeds, the local variables $w_i$ will get closer to each other and will track increasingly better the gradient of the overall cost function $f(x)$.

\renewcommand{\thealgorithm}{}
\begin{algorithm}[b]
\caption{Network-GIANT}\label{alg_NetworkGIANT}
\begin{algorithmic}
\Initialize{
     Arbitrary $x_i^0 \in \mathbb{R}^d$, \\
     $w_i^0 = \nabla f_i (x_i^{-1}) = \mathbb{0} \quad \forall i \in \mathcal{N}$, \\
    $\epsilon > 0$, consensus matrix $P$.}
\For{each iteration $k = 0, 1, 2, \dots$}
    \For{each agent $i = 1, \dots, n$}
    
        \State $v^{k+1}_i = w_i^k + \nabla f_i (x_i^k) - \nabla f_i (x_i^{k-1})$
        \State $w^{k+1}_i = \sum_{j=1}^n p_{ij} v^{k+1}_j$
        \State $u_i^{k+1} =  x_i^k - \epsilon \nabla^2 f_i (x_i^{k})^{-1} w_i^{k+1}$
        \State $x^{k+1}_i = \sum_{j=1}^n p_{ij} u^{k+1}_j$     
    \EndFor
\EndFor
\end{algorithmic}
\end{algorithm}

With these ingredients at hand, we are ready to state our distributed second-order algorithm for general networks, which we refer to as Network-GIANT. The local variables of the agents are initialized with arbitrary values of  $x_i^0 \in \mathbb{R}^d$ and $w_i^0 = \nabla f_i (x_i^{-1})$ $\forall i \in \mathcal{N}$. At each iteration $k$, each agent $i$ computes the gradient $\nabla f_i (x_i^k)$ of its local function and updates $w_i^k$ using \eqref{eq:w_gradient_tracking}, which involves a consensus step. The resulting $w_i^{k+1}$ is combined with the local Hessian to obtain an approximate Newton direction and descend along it, obtaining
\begin{equation}
    u_i^{k+1} =  x_i^k - \epsilon \nabla^2 f_i (x_i^{k})^{-1} w_i^{k+1} .
    \label{eq:u_update_rule}
\end{equation}

The iteration terminates with a consensus step on the direction $u_i^{k+1}$, whose result provides the value of $x^{k+1}_i$. Compared with GIANT, which employs a line search procedure, in \eqref{eq:u_update_rule} we use a fixed step size $\epsilon > 0$. This implies that at each iteration to only two $d$-dimensional vectors per agent are transmitted, resulting in an appealing communication cost.

\begin{remark}{}
Aligning with similar literature, e.g. ESOM-$K$ \cite{7576649}, Network-DANE for different values of $K$ \cite{li2020communication}, NN-$K$ \cite{mokhtari2016network}, also in Network-GIANT it is possible to run multiple consensus steps within each consensus block. For example, the single consensus step in the last line of the pseudocode can be replaced with $K \geq 1$ consecutive consensus rounds. This is equivalent to a single consensus step in which the mixing matrix $P^K$ is used. The parameter $K$ adds flexibility to the algorithm, introducing a trade-off between communication cost, wall-clock time and number of iterations needed for convergence.
\end{remark}


\section{Convergence Analysis}
In this section we prove the semi-global exponential convergence of Network-GIANT to the exact solution of the optimization problem \eqref{eq:problem_formulation}. The following assumption formalizes the convexity and regularity properties of the objective functions.

\begin{assumption} [Strong Convexity and Smoothness]
Let the local costs $f_i(x)$ $\forall i \in \mathcal{N}$ be two times differentiable, $\mu$-strongly convex and with $L$-Lipschitz and continuously differentiable gradients. This assumption, which is standard in the field of convex optimization, tells that there exist positive constants $\mu, L$ such that $\mu I \leq \nabla^2 f_i (x) \leq L I$, $\forall x \in \mathbb{R}^d$. To simplify the proof of Theorem \ref{theorem:convergence_Network_GIANT} we further characterize the objective functions assuming that the inverse Hessian matrices are continuously differentiable functions.
\label{assumpt:convexity_smoothness}
\end{assumption}

Our proof is based on the principle of time-scales separation for discrete-time systems \cite{ZO-JADE}. To keep the manuscript self-contained and for a better understanding of the proof of our main result, we provide a statement of such theorem.

\begin{theorem} [Separation of time-scales \cite{ZO-JADE}]
The origin of the autonomous discrete-time system
\begin{equation*}
\left\{\begin{array}{l}
z^k=z^{k-1}+\epsilon \phi \left( z^{k-1}, y^{k-1} \right) \\
y^k=\varphi \left(y^{k-1}, z^{k-1} \right) + \epsilon \Psi \left(y^{k-1}, z^{k-1} \right)
\end{array}\right.
\end{equation*}
is semi-globally exponentially stable for a suitable choice of the parameter $\epsilon > 0$, provided that all the following conditions are verified:
\textit{(a)} The functions $\varphi$, $\Psi$ and $\phi$ are locally uniformly Lipschitz.
\textit{(b)} There exists $y^\star$ such that $y^\star(z)=\varphi \left(y^\star(z), z \right)$ for all $z$.
\textit{(c)} At the origin, $\phi\left(0, y^{*}(0)\right)=0$ and  $\Psi\left(y^{*}(0), 0\right)=0$.
\textit{(d)} There exists a twice differentiable function $V (z)$ and positive constants $c_1,c_2,c_3,c_4$ such that
    \begin{equation*}
    \begin{gathered}
    c_{1}\|z\|^{2} \leq V(z) \leq c_{2}\|z\|^{2}, \\
    \frac{\partial V}{\partial z} \phi\left(z, y^{*}(z)\right) \leq-c_{3}\|z\|^{2},
    \qquad
    \left\|\frac{\partial V}{\partial z}\right\| \leq c_{4}\|z\| .
    \end{gathered}
    \end{equation*}
\label{theorem:sep_time_scales}
\end{theorem}

The characterization ”semi-globally exponentially stable” means that Network-GIANT exhibits a linear convergence rate when initialized within a neighborhood of the optimal solution. A more rigorous definition of can be found in Proposition 8.1 in \cite{bof2018lyapunov}.
Below we present our main result, which provides theoretical guarantees on the stability and performance of the proposed algorithm. Since the parameter $K$ does not affect the convergence proof, for simplicity we set $K=1$.

\begin{theorem} [Convergence of Network-GIANT]
Under Assumption \ref{assumpt:convexity_smoothness}, there exists a positive $\bar{\epsilon}$ such that for any $\epsilon \leq \bar{\epsilon}$ Network-GIANT converges semi-globally and exponentially fast to the exact solution of \eqref{eq:problem_formulation}.
\label{theorem:convergence_Network_GIANT}
\end{theorem}

\begin{proof}
The proof follows the same steps of Theorem 2 in \cite{ZO-JADE}, but presents different challenges. Our goal is to rewrite Network-GIANT as a discrete-time system compatible with Theorem \ref{theorem:sep_time_scales}, and then show that all the assumptions of the latter are satisfied. Without loss of generality we translate the objective function by an offset, so that the solution $x^\star = \text{argmin} f(x) = \mathbb{0}$ and $f(x^\star) = \mathbb{0}$.
To obtain a compact notation we define the stacks of variables $\in \mathbb{R}^{n \times d}$
\begin{equation*}
\begin{gathered}
x^k \coloneqq [ x_1^k, \dots, x_n^k ]^T,
\quad
\nabla f(x^k) \coloneqq [ \nabla f(x_1^k), \dots, \nabla f(x_n^k) ]^T,
\\
\nabla^2 f(x^k)^{-1} w^{k+1} \coloneqq [\nabla^2 f_1(x_1^k)^{-1} w_1^{k+1}, \dots]^T .
\end{gathered}
\end{equation*}
We decompose $x(t)$ into the sum of its mean $\bar{x} = \frac{1}{n}  \mathbb{1} \mathbb{1}^T x$, which has all equal rows, and $\tilde{x} = x - \bar{x}$.
Using the doubly stochasticity of the matrix $P$, it is easy to see that Network-GIANT is equivalent to the discrete-time dynamical system identified by the slow dynamics
\begin{equation*}
\begin{aligned}
\bar{x}^{k+1} &= \bar{x}^k + \epsilon \left( - \frac{1}{n}  \mathbb{1} \mathbb{1}^T  \nabla^2 f(\bar{x}^k + \tilde{x}^k)^{-1} w^{k+1} \right) \\
&= \bar{x}^k + \epsilon \phi(\bar{x}^k, \xi^k)
\end{aligned}
\label{eq:dyn_sys_slow}
\end{equation*}
and the fast dynamics
\begin{equation*}
\begin{aligned}
g^{k+1} &= \varphi_g = \nabla f(\bar{x}^k + \tilde{x}^k), \\
w^{k+1} &= \varphi_w = P \left( w^k + \nabla f(\bar{x}^k + \tilde{x}^k) - g^k \right), \\
\tilde{x}^{k+1} &= \left( I - \frac{1}{n} \mathbb{1} \mathbb{1}^T \right) P \left( x^k - \epsilon \nabla^2 f(x^k)^{-1} w^{k+1} \right) \\
&= \varphi_{\tilde{x}} + \epsilon \Psi(\bar{x}^k, \xi^k), \\
\text{where} &
\quad
\xi = \left\{ g, w, \tilde{x} \right\},
\quad
\varphi_{\tilde{x}} = \left( P - \frac{1}{n} \mathbb{1} \mathbb{1}^T \right)\tilde{x}^k, \\
\Psi(\bar{x}^k, \xi^k) &= \left(\frac{1}{n} \mathbb{1} \mathbb{1}^T - P \right) \nabla^2 f(\bar{x}^k + \tilde{x}^k)^{-1} w^{k+1}.
\end{aligned}
\end{equation*}

Using the definition of $w^{k+1}$ we note that the system is autonomous, and setting $z = \bar{x}$, $y = \xi$ and $\varphi = \varphi_g + \varphi_w + \varphi_{\tilde{x}}$ we get the same decoupled structure of the system considered in Theorem \ref{theorem:sep_time_scales}. Below we show that all the conditions of such Theorem are verified.

Since by $w$ is a linear combination of continuously differentiable and Lipschitz continuous terms, we have that $w$ is itself continuously differentiable and Lipschitz.
Consequently, also the product of continuously differentiable functions $\nabla^2 f(x^k)^{-1} w^{k+1}$ is continuously differentiable and locally Lipschitz. The above facts guarantee that the functions $\varphi$, $\phi$ and $\Psi$ are continuously differentiable and Lipschitz.

The boundary-layer system $\xi^k=\varphi \left(\xi^{k-1}, \bar{x} \right)$, obtained by setting $\epsilon = 0$, reaches the steady-state value $\xi^\star(\bar{x}) \coloneqq \{ g^\star = w^\star = \nabla f(\bar{x}),\ \tilde{x}^\star = \mathbb{0} \}$ for any choice of $\bar{x}$.
Recalling that $x^\star = \mathbb{0}$, by the first-order optimality condition $\nabla f(\mathbb{0}) = \mathbb{0}$, which implies $\phi \left(\mathbb{0}, \xi^\star(\mathbb{0}) \right) = \mathbb{0}$ and $\Psi\left(y^{*}(\mathbb{0}), \mathbb{0}\right)=\mathbb{0}$. 
Finally, a suitable choice for the Lyapunov function is $V(\bar{x})= f(\bar{x})$. Indeed, using Taylor expansion, the fact that $f(\mathbb{0}) = \mathbb{0}$ and Assumption \ref{assumpt:convexity_smoothness} the following inequalities hold.
\begin{equation*}
\frac{\mu^2}{2} \left\| \bar{x} \right\|^2 \leq
V(\bar{x}) \leq \frac{L^2}{2} \left\| \bar{x} \right\|^2,
\quad
\left\| \frac{\partial V}{\partial \bar{x}} \right\|
\leq L \left\| \bar{x} \right\|,
\end{equation*}
\begin{equation*}
\begin{aligned}
\frac{\partial V}{\partial \bar{x}} \phi \left(\bar{x}, \xi^\star(\bar{x}) \right)
&= - \nabla f(\bar{x})^T \left( \frac{1}{n} \sum_{i=1}^n \left[ \nabla^2 f_i(\bar{x})^{-1} \right] \right) \nabla f(\bar{x}) \\
&\leq - \frac{1}{L} \left\| \nabla f(\bar{x}) \right\|^2
\leq - \frac{\mu^2}{L} \left\| \bar{x} \right\|^2.
\end{aligned}
\end{equation*}
Since all the required conditions are satisfied, we can invoke Theorem \ref{theorem:sep_time_scales} to complete the proof.

\end{proof}

\begin{remark}
    The proof of Theorem \ref{theorem:convergence_Network_GIANT} utilizes the time-scales separation principle which allows to prove linear convergence and get concise proofs. However, this principle is primarily an existential theorem and makes it difficult to determine the exact convergence rate coefficient. While we do not explicitly compare our theoretical convergence rate with other algorithms, it is important to note that such rate can be at most linear for all the consensus-based algorithms due to the bottleneck imposed by the mixing matrices.
\end{remark}

\section{Numerical Analysis}
\begin{figure*}[t!]
    \centering  
    
    \begin{subfigure}[b]{\textwidth}
        \centering
        \includegraphics[width=0.95\textwidth]{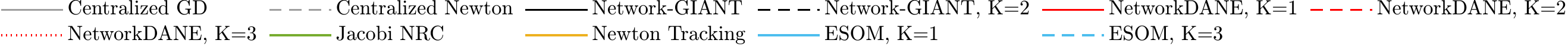}
        \vspace{0.3cm}
    \end{subfigure}
          
     \begin{subfigure}[b]{\textwidth}
        \centering
        \includegraphics[width=0.24\textwidth]{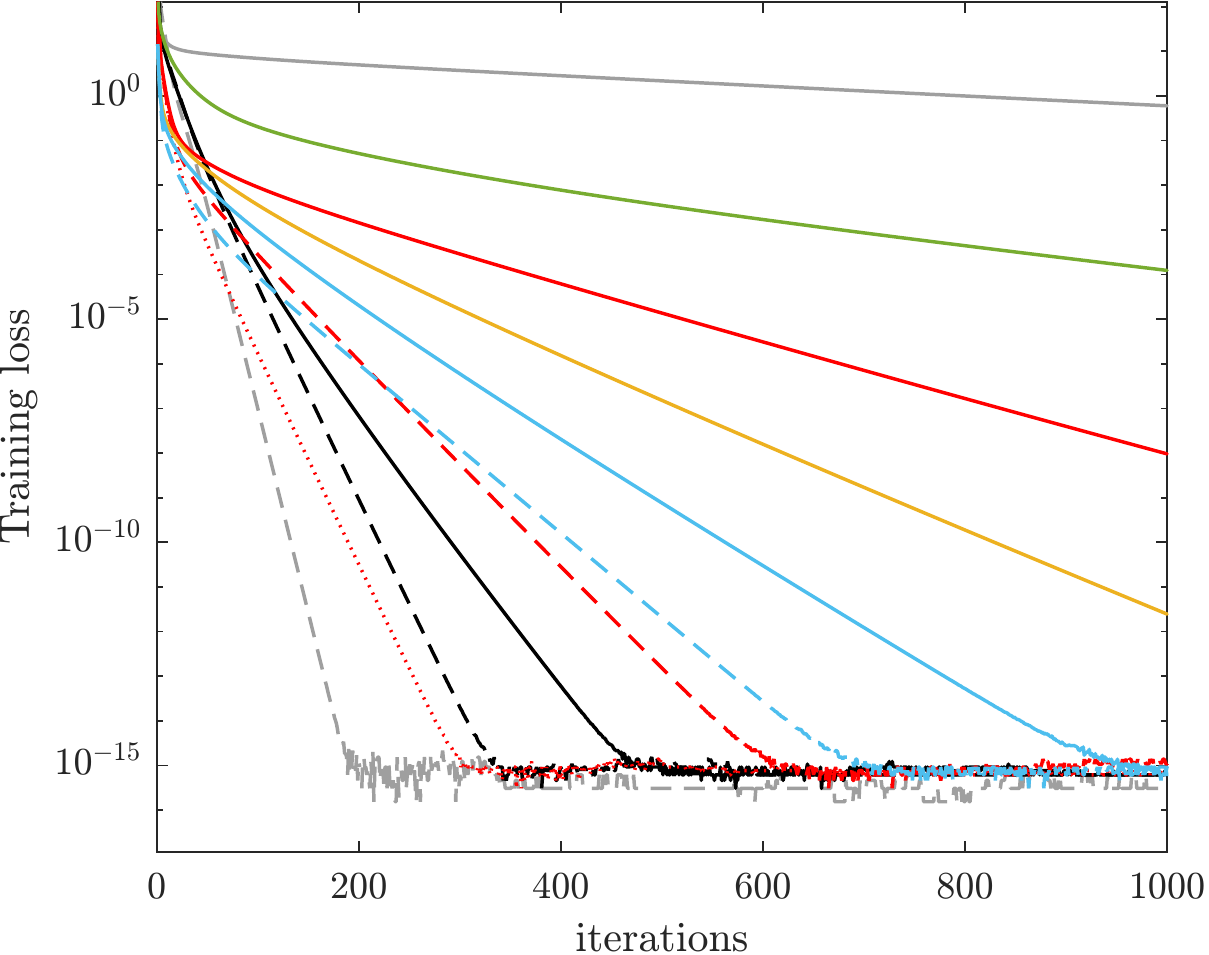}
        \hfill
        \includegraphics[width=0.24\textwidth]{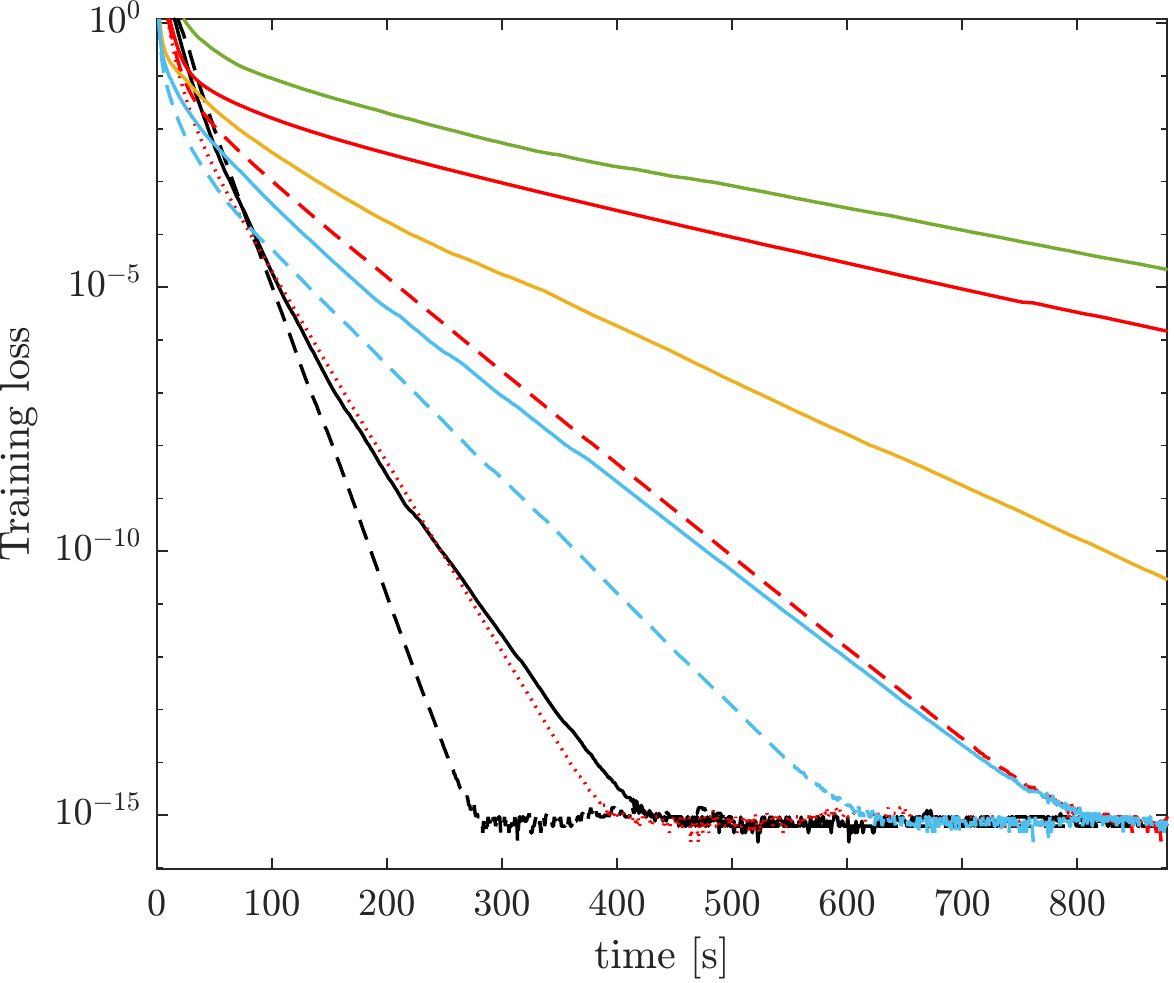}
        \hfill
        \includegraphics[width=0.24\textwidth]{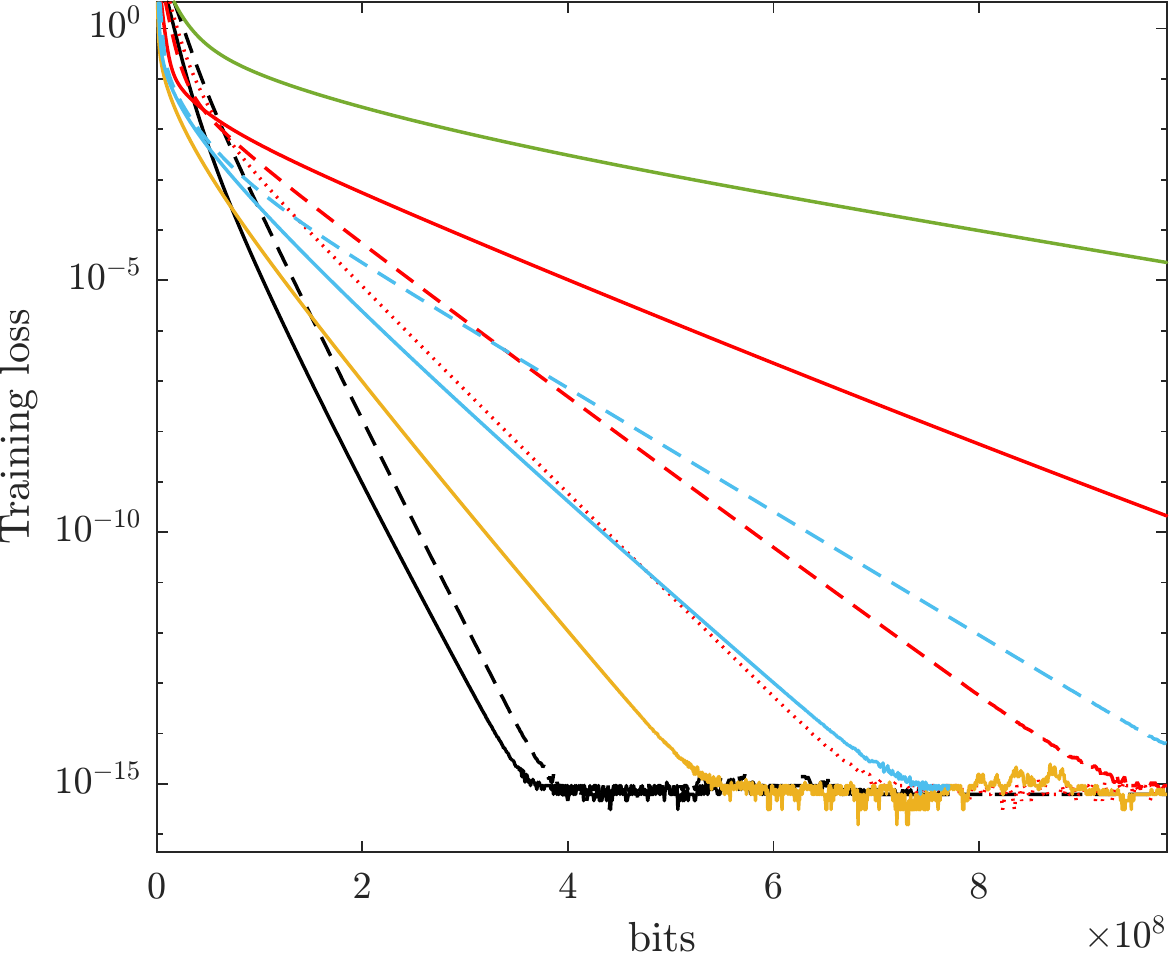}
        \hfill
        \includegraphics[width=0.24\textwidth]{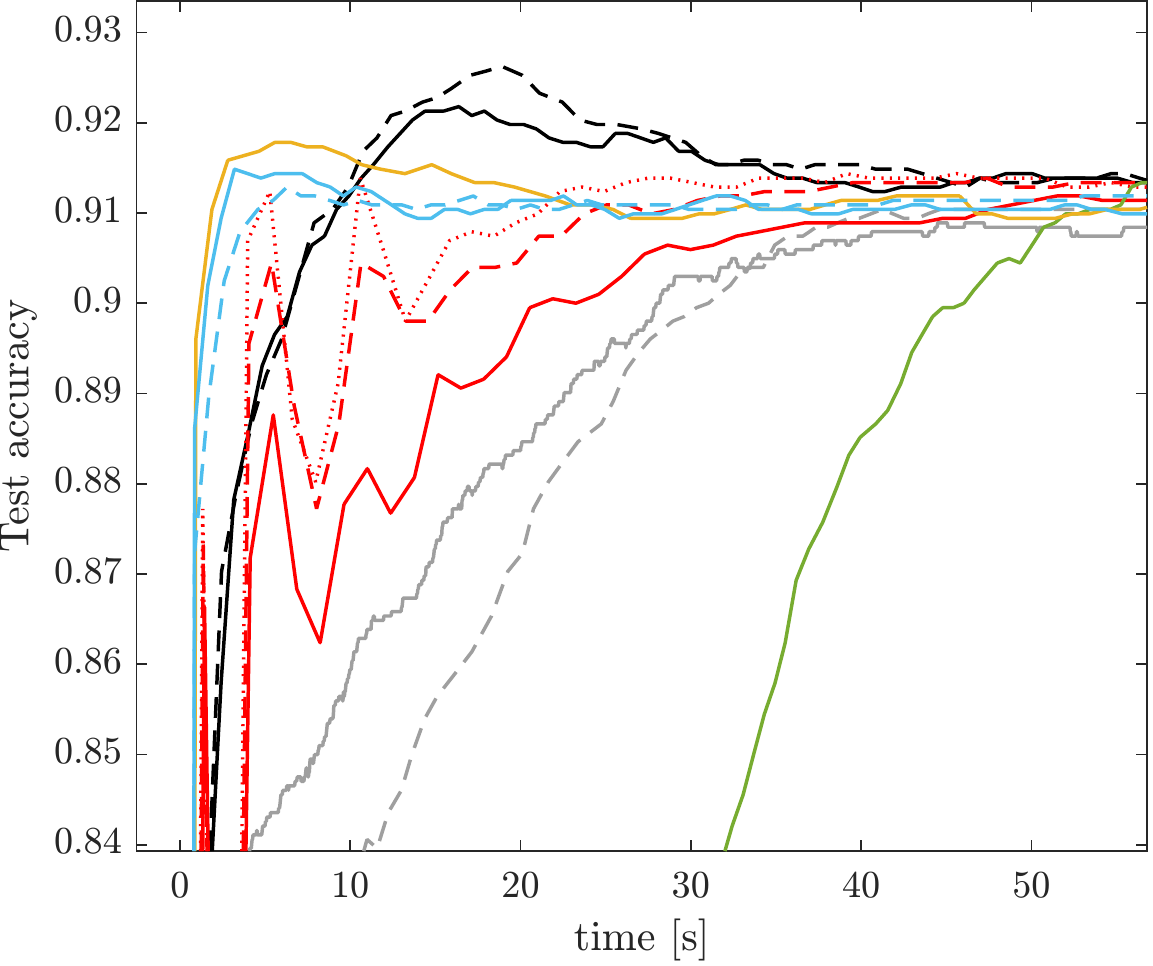}
        \caption*{(a) Results on the dataset MNIST.}
     \end{subfigure}

 \begin{subfigure}[b]{\textwidth}
        \centering
        \includegraphics[width=0.24\textwidth]{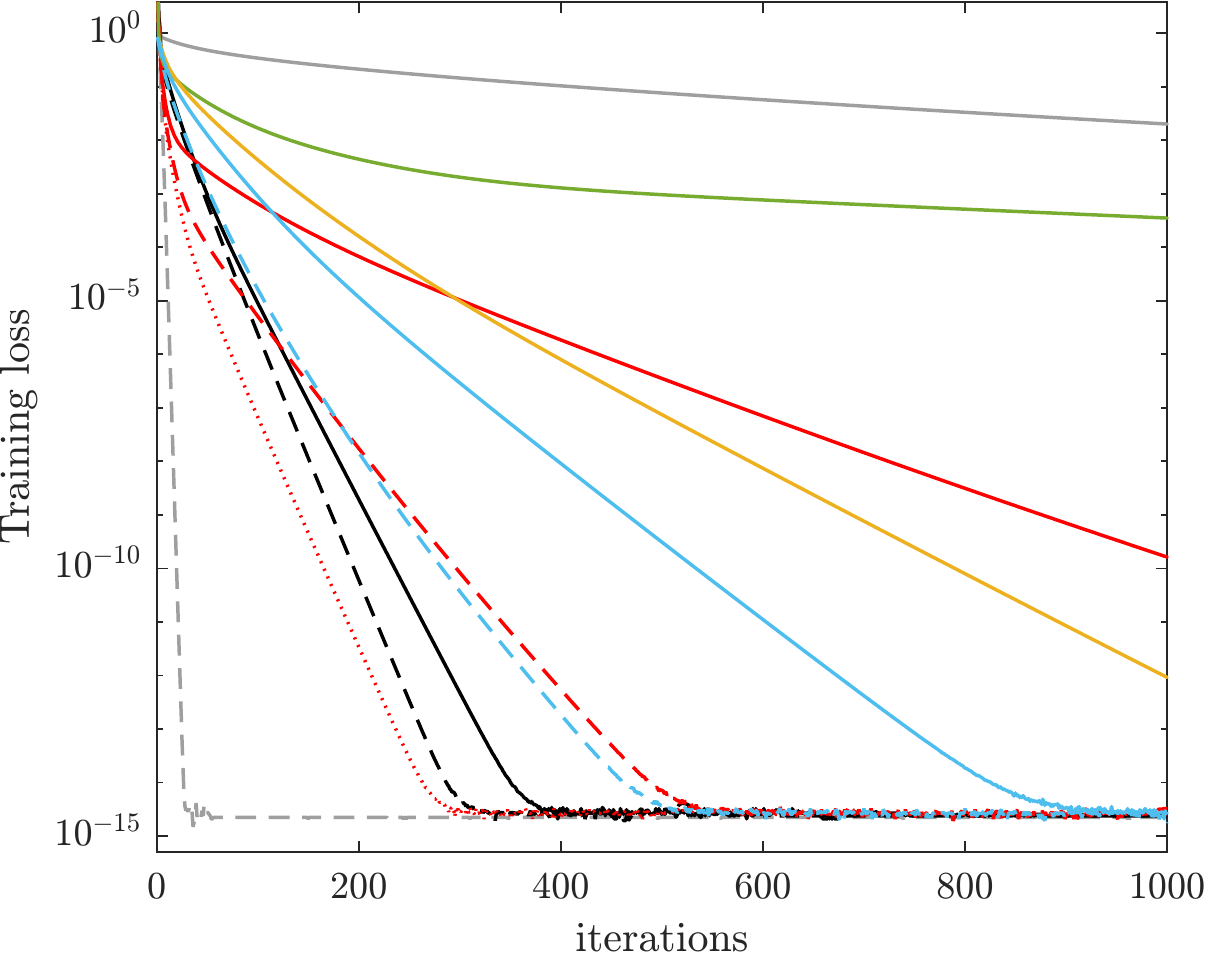}
        \hfill
        \includegraphics[width=0.24\textwidth]{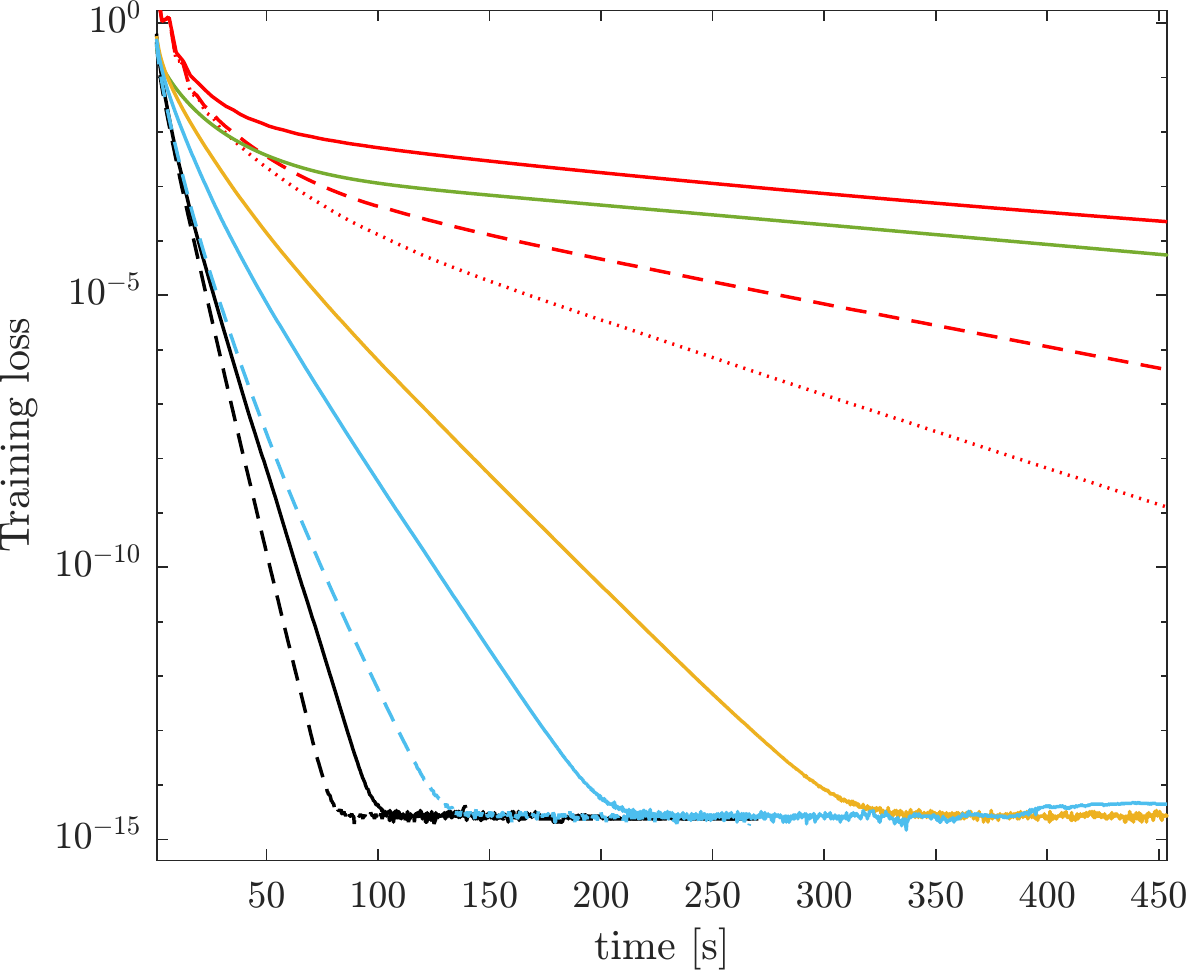}
        \hfill
        \includegraphics[width=0.24\textwidth]{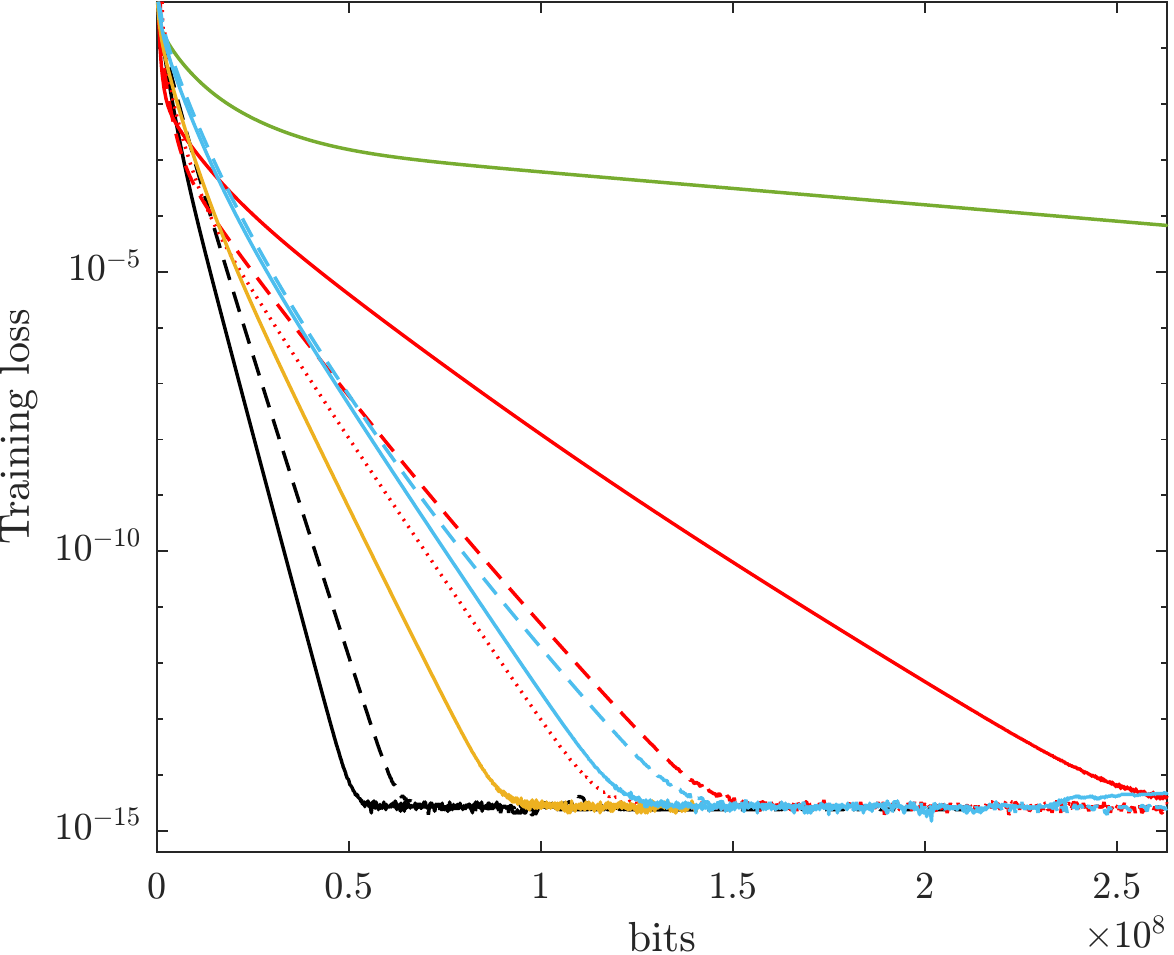}
        \hfill
        \includegraphics[width=0.24\textwidth]{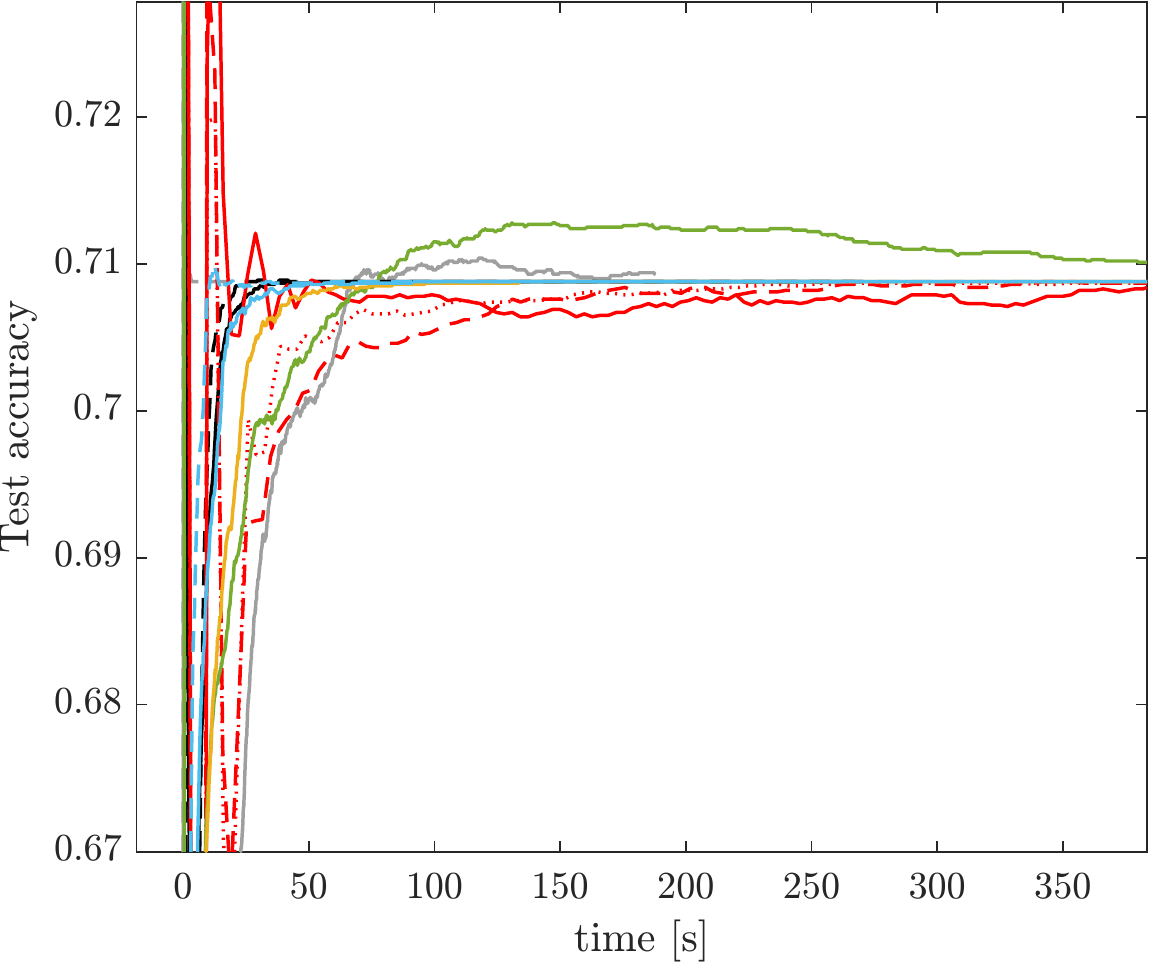}
     \caption*{(b) Results on the dataset Covertype.}
     \end{subfigure}

     \caption{In both (a) and (b) the first three figures from the left show the performance on the training set, while the last figure on the right displays the transient evolution of the accuracy on the test set. The convergence is plotted against multiple metrics, namely the number of iterations, the run-time and the overall communication cost in bits}

    \label{fig}  
\end{figure*}

In this section we practically assess the performance of Network-GIANT through numerical experiments and compare it with several other state-of-the-art distributed algorithms. In particular, we perform distributed binary classification between two classes of the well-known datasets MNIST \cite{lecun2010mnist} and Covertype \cite{dheeru2017uci} via logistic regression. For the MNIST dataset we use Principal Component Analysis to reduce the feature size of to $d=300$, while for the Covertype dataset we keep the original dimension $d=54$.
We consider a network of $20$ nodes and we build the mixing matrix $P$ using the Metropolis-Hastings weights \cite{xiao2007distributed}.
Each node $i$ is given $m_i$ data samples $(c^{(i)}, y^{(i)}) \in \mathbb{R}^d \times \{-1,+1\}$, where the first element in the tuple is the predictor and the second is the target response. The local cost functions are the regularized log-losses
\[
f_i(x) = \frac{1}{m_i}\sum_{j=1}^{m_i} \log(1+e^{-y^{(i)}_j(x^T c^{(i)}_j)}) + \frac{\lambda}{2} \|x\|^2,
\]
while the training loss shown in the plots is
\begin{equation*}    
\frac{\frac{1}{n}\sum_{i=1}^{n}f_i(x)-f(x^*)}{f(x^*)} .
\end{equation*}

We test Network-GIANT both in its vanilla version and when $K=2$ consecutive consensus steps on $\{u_i\}_{i=1}^n$ are performed.
We compare them with the following fully distributed algorithms: Network-DANE \cite{li2020communication} and ESOM-$K$ \cite{7576649}, both with different values $K$ of mixing rounds; NRC \cite{varagnolo2015newton} in its Jacobi version, which keeps the communication cost linear in $d$; Newton Tracking \cite{zhang2021newton}. We do not compare Network-GIANT to any federated learning algorithm, not even GIANT, because the results would depend on the chosen network topology. Also, we do not consider the inexact algorithms Network Newton \cite{mokhtari2016network} and DQN \cite{bajovic2017newton} which only converge to a neighborhood of the solution. Rather, we include in the comparisons the standard centralized gradient descent (GD) with backtracking line search and the centralized Newton-Raphson method, to be used as a reference.
For all the algorithms we tune the hyper-parameters as suggested in the respective papers, selecting the values that provide the best performance.

Fig. \ref{fig} provides interesting insights into the strengths and weaknesses of the algorithms under examination, allowing a complete performance comparison that takes into account all the main evaluation metrics. When the convergence rate is expressed in terms of either run-time or communication cost, Network-GIANT is superior to all its competitors. This consistency does not hold for the other algorithms: for example, ESOM outperforms Newton Tracking in terms of number of iterations and run-time, but the opposite is true if we look at the communication cost.
When looking at the number of iterations, Network-DANE with $K=3$ converges slightly faster than the others. However, each iteration of Network-DANE involves an inner minimization problem which in general is more computationally demanding compared to the other algorithms, even for $d$ in the size of hundreds. Finally, the relatively poor performance of Jacobi NRC can be attributed to the fact that this algorithm exploits only the diagonal of local Hessian matrices, leading to slow convergence for skewed objective functions.

The numerical experiments also reveal the importance the parameter $K$, which dramatically affects the ranking of the various algorithms. In particular, for Network-DANE bigger values of $K$ always seem beneficial with respect to all metrics, up to a certain threshold. Instead, for Network-GIANT and ESOM-$K$ we observe a trade-off: higher values of $K$ decrease both the number of iterations and the run-time needed for convergence, but at the same time increase the total number of bits transmitted between the nodes.

For completeness we also show the performance on the test set, recalling that generalization to out-of-sample data is a bigger problem outside the scope of this work. All the algorithms, with the exception of Jacobi NRC, achieve similar accuracy at a comparable rate, exhibiting some oscillations that may be due to the statistical difference between training and test set.

\section{Conclusions}
This paper presents the development of a highly efficient network learning algorithm, named Network-GIANT. The proposed algorithm eliminates the need for a parameter server by applying the GIANT algorithm to a network learning framework. It leverages the concepts of average consensus and gradient tracking to achieve an approximate Newton-type peer-to-peer distributed optimization. The convergence analysis of Network-GIANT indicates a semi-global exponential convergence to the exact solution. The paper provides an empirical comparison of Network-GIANT with other state-of-the-art network learning algorithms, such as Network-DANE, Jacobi NRC, Newton Tracking, and ESOM. The results demonstrate that Network-GIANT outperforms the other network algorithms in terms of convergence, communication cost and computation performance. We foresee interesting possibilities for future work, such as the analysis of robustness to disturbances and noise injection, and the introduction of quantization to further reduce communication costs.

\bibliographystyle{IEEEtran}
\bibliography{utils/refs}

\end{document}